\documentclass[11pt,a4paper]{amsart}

\usepackage{amsmath,amsfonts,amsthm,amssymb}

\title{On Fano indices of $\Q$-Fano 3-folds}
\author{Kaori Suzuki}
\date{}
\address{Department of Mathematical Sciences, University of
Tokyo, 
Komaba, Meguro, Tokyo 153--8914, Japan}
\email{suzuki@ms.u-tokyo.ac.jp}

\theoremstyle{definition}
\newtheorem{definition}{Definition}[section]

\theoremstyle{plain}
\newtheorem{thm}[definition]{Theorem}
\newtheorem{prop}[definition]{Proposition}
\newtheorem{lemma}[definition]{Lemma}
\newtheorem{cor}[definition]{Corollary}

\newtheorem{cl}[definition]{Claim}

\theoremstyle{remark}
\newtheorem{rem}[definition]{Remark}

 \newcommand{\C}{\mathbb C}
 \newcommand{\Q}{\mathbb Q}
 
 \newcommand{\Z}{\mathbb Z}
 \newcommand{\E}{\mathcal E}
 \newcommand{\F}{\mathcal F}

\DeclareMathOperator{\lcm}{lcm}

 \def\coh2#1{{\mathit{H}^2({#1},\Z)}{}}
\numberwithin{equation}{section}

\newcommand{\PP}{\mathbb P}
\newcommand{\Oh}{\mathcal O}

\begin{document}

\begin{abstract} 
We shall give the best possible upper bound of 
the Fano indices together with a characterization of those $\Q$-Fano 
3-folds which attain the maximum in terms of graded rings. 
\end{abstract}

\maketitle

\renewcommand{\leftmark}{K. Suzuki}
\renewcommand{\rightmark}{Fano 3-folds}
\setcounter{section}{-1}

\section{Introduction}


$\Q$-Fano 3-folds play important roles in birational algebraic 
geometry. They have been studied by several authors since G. Fano. 
In this paper, 
we study $\Q$-Fano 3-folds from the view of their 
Fano indices (See definition 0.2 below) and graded rings. More concretely, 
we give an optimal upper bound for the Fano indices and also characterize 
those $\Q$-Fano 
3-folds which attain the maximum in terms of graded rings (Theorem 0.3). 
Throughout this paper, 
we work over the complex number field $\C$. 

\par
\vskip 4pt 

\begin{definition}
Let $X$ be a normal projective 3-fold. We call $X$ a $\Q$-Fano 3-fold if: 
\begin{enumerate}
\item[(1)] $X$ has only $\Q$-factorial terminal singularities;
\item[(2)] the anti-canonical (Weil) divisor $-K_{X}$ is ample; and
\item[(3)] $\rho(X) = 1$, $\rho(X)$ is the Picard number of $X$. 
\end{enumerate}
\end{definition} 

Let $X$ be a $\Q$-Fano 3-fold. There are two important indices of $X$: 

\begin{definition} 

 We define the Gorenstein index $r=r(X)$ and the Fano inex $f=f(X)$ by  
\[
r(X) := \text{min}\, \{n \in \Z_{> 0}\, 
\vert nK_{X}\, \text{is Cartier}\}\, ;
\]
 
\[
f(X) := \max\, \{m \in \Z_{> 0}\, 
\vert K_{X} = mA\,~\textrm{for some integral Weil divisor}\,~A\}\, .
\]
Here the equality $K_{X} = mA$ means that $K_X-mA$ is linear equivalent to $0$. 
If $-K_{X} = f(X)A$, we call $A = A_{X}$ a primitive 
Weil divisor.

\end{definition} 

In earlier works of Shokurov,
Alexeev, Iskovskikh, Prokhorov, Sano, Mella and others the Fano
index was defined in a different way, as the maximal rational such
that $-K_X \equiv rH$ for some ample Cartier divisor $H$. Note 
that our definition is different
from the one used by previous authors.

Although the Gorenstein indices do not appear in the statement of main results, 
they play crucial roles in the proof (See section 2). 
\par
\vskip 4pt
Our main result is as follows: 

\begin{thm}
Set 
\[
\F := \{n \in \Z_{>0} \vert 1 \leq n \leq 11,~\textrm{or}~ 13, 17, 19\} =  \{1, 2, \cdots , 9, 10, 11, 13, 17, 19\}\,\, .
\] 
Let $X$ be a $\Q$-Fano $3$-fold of Fano index $f(X)$ and $A = A_{X}$ 
 a primitive Weil divisor. Then: 
\begin{enumerate} 
\item[(1)] $f(X) \in \F$. In particular, $f(X) \leq 19$. 
\item[(2)] If $X$ is the $\PP(3, 4, 5, 7)$, $f$ attains the maximum 19. 

In addition, for any $X$ with $f(X) = 19$, 
the Hilbert series of $(X, A_{X})$ coincides with that of 
$(\PP (3, 4, 5, 7), \Oh(1))$, i.e. 
\[
\sum_{n \geq 0} h^{0}(X, \Oh_{X}(nA_{X}))t^{n} = 
\frac{1}{(1-t^{3})(1-t^{4})(1-t^{5})(1-t^{7})}\, .
\] 
\item[(3)] Each element of $\F$ except possibly $10$ is realized 
as a Fano index. 
\end{enumerate}
\end{thm} 
We expect the uniqueness of 
$X$ with $f(X) = 19$, to which the second statement of (2) provides a 
supporting 
evidence. For the statement (3), we shall 
construct desired examples as hypersurfaces in suitable 
weighted projective spaces (Section 2). As a by-product, it turns out that 
for each 
$f \in \F - \{10\}$, there is a $\Q$-Fano 3-fold $X$ of Fano index 
$f$ with only cyclic quotient terminal singularities. We expect that there is no 
Fano 3-fold of Fano index $10$. We hope to return back this problem in future. 
\par
\vskip 4pt
Our proof is based on an (effective version of) 
boundedness theorem of $\Q$-Fano 3-folds due to Kawamata [Ka3] 
(Theorem 1.7), the singular Riemann-Roch formula by Reid [Re] (Theorem 1.4). 
In order to make our estimate optimal, we also use 
computer programs called Magma [Ma] 
at the final stage. We emphasize that our use of computer programs involves 
nothing more 
than addition, subtraction, multiplication and division of reasonable amount 
of positive integers, which, in principle, can be done also by hand. We 
collect the necessary programs in the appendix for interesting 
readers. 
\par
\vskip 4pt

\subsection*{Acknowledgement}
The author would like to express her deep gratitude to Professor Miles Reid for his valuable comments and unceasing encouragement. 
She thanks Professors Toshiyuki Katsura and Youichi Miyaoka and Doctor Hiromichi Takagi for helpful discussions. 
The Magma programs exploited in this paper were supported by Doctor Gavin Brown at University Warwick. Special thanks go to him 
for teaching how to use the programs. 
The author is grateful to Professor Sinobu Hosono for verifying some cases of 
our results on Mathematica by translating our computer programs. 
The author gratefully thanks to her PhD supervisor Professor Keiji Oguiso for his warm encouragement. She also would like to thank the referee, who read 
the manuscript carefully and gave useful suggestions. This research was 
partially 
supported by the 21st Century COE Program 
at Graduate School of Mathematical Sciences, the University of Tokyo.

\section{Preliminaries}

In this section, we recall the notion called the basket of 
singularities after Mori [Mo] and Reid [Re], and 
two fundamental theorems, namely, 
the singular Riemann-Roch theorem due to Reid [Re] 
and the boundedness theorem of $\Q$-Fano 3-folds due to Kawamata [Ka3]. 
These two theorems will be essential for our study. 
\par 
\vskip 4pt

Let $(U, P)$ be a germ of a 3-dimensional terminal singularity of index 
$r = r_{P} > 1$. It is known by [Mo] that if $(U, P)$ is not a quotient 
singular point, $(U, P)$ can be deformed to a (unique) collection of a 
finite number of terminal quotient singularities, say 
$\{P_{k}\}_{k=1}^{n_{P}}$. 
Write a type of singularity of $P_{k}$ as 
\[
\frac{1}{r_{P, k}}(1, a_{P, k}, -a_{P, k})\,~\textrm{or}~\, [r_{P, k}, a_{P, k}]\, .
\]
Then $(r_{P, k}, a_{P, k}) = 1$ and $r_{k} \geq 2$, $n_{P} \geq 2$ and 
$r_{P} = \,\lcm\, (r_{P, k})$ hold\footnote{By [Mo], it is known that 
except one exceptional case where 
$r_{P} = r_{P, 1} = 4$ and $r_{P, i} = 2$ ($i \geq 2$), one has also 
$r_{P, k} = r_{P}$ for all $k$.}. We call the set 
$\mathcal{B}(U, P) := \{P_{k}\}_{k=1}^{n_{P}}$ the basket of singularities of 
$(U, P)$. When $(U, P)$ is already a quotient singular point, we regard 
$\{(U, P)\}$ itself as the basket of singularities of $(U, P)$. Since the $3$-dimensional 
terminal 
singularities are isolated singularities, we can speak of the basket of 
singularities in the global case: 

\begin{definition}
Let $X$ be a terminal 3-fold. Let $\{P_{i}\}_{i=1}^{m}$ be the set of 
singular points of $X$ with $r(P_{i}) \geq 2$ and $P_{i} \in U_{i} (\subset X)$ be a 
small analytic neighborhood of $P_{i}$. 
Then, we call the disjoint union 
$\cup_{i=1}^{m} \mathcal{B}(U_{i}, P_{i})$ the basket of singularities of $X$. 
\end{definition} 

We often describe the basket of singularities of $X$ by listing up the type of each 
point in the basket, like 
\[
\mathcal{B}(X) = 
\{[2, 1], [2, 1], [2, 1], [5, 2], [7, 2]\}\, .
\]

\par
\vskip 4pt
As it will be reviewed below, several important invariants of $\Q$-Fano 
3-folds depend only on the basket of singularities (but not on the actual set 
of singularities of $X$). However, we should also notice that the basket of 
singularities encodes {\it no information} about Gorenstein singular points. 
For instance, if 
\[
\mathcal{B}(X) = \{[3, 1], [4, 1], [5, 2], [7, 2]\}\, , 
\] 
then one can deduce that $X$ has four quotient 
singular points of indicate types by using the facts explained above and 
in the footnote. 
However one can not say anything about 
Gorenstein singular points of $X$. 
\par
\vskip 4pt 
Let $X$ be a $\Q$-factorial terminal 3-fold and $P \in X$ be a singular 
point of local Gorenstein index $r_{P} > 1$. Here the local Gorenstein index $r_{P}$ is defined to be 
the smallest positive integer $r_{P}$ such that 
$r_{P}K_{X}$ is Cartier at $P$. 
In particular, by [Ka2, Corollary 5.2], the local class group 
at $P$ is isomorphic to the cyclic group $C_{r_{P}}  := \Z /r_{P}$ and $K_{X}$ is a 
generator of the local class group. Then, for each Weil divisor $D$ on $X$, 
there is a unique integer 
$i:= i(P, D) \in [0, r_{P}-1]$ such that $D = iK_{X}$. We call the integer 
$i = i(P, D)$ the local index of $D$ at $P$. When $D$ can be deformed 
under the general deformation $\mathcal{U} \rightarrow \Delta$ of the germ 
$P \in U$, we can write a similar equation in terms of the basket of 
singularities $\{P_{k} = [r_{k}, a_{k}]\}$ of $(U, P)$, i.e. there is a 
unique integer $i := i(P_{k}, D_{k}) \in [0, r_{k} -1]$ such that 
$D_{t} = iK_{\mathcal{U}_{t}}$ around $P_{k} \in \mathcal{U}_{t}$. By abuse of notation, 
we also call the integer $i := i(P_{k}, D_{k})$ the local index of $D$ at 
(the basket point) $P_{k}$.

\begin{lemma} 
Let $X$ be a $\Q$-Fano $3$-fold. Let $\mathcal{B}(X) = \{[r_{k}, a_{k}]\}_{k=1}^{m}$ be the basket of singularities of $X$, $r := r(X)$ be 
the Gorenstein index of $X$, 
$f= f(X)$ be the Fano index and $A$ be a primitive Weil divisor. Then 
\begin{enumerate} 
\item[(1)] $r$ is the least common multiple of $\{r_{k}\}_{k=1}^{m}$, i.e. 
$r = \lcm\, (r_{k})_{k=1}^{m}$. 
\item[(2)] $A^{3} = n/r$ for some positive integer $n$. 
\item[(3)] $r_{k}$ and $f$ are coprime, i.e. $(r_{k}, f) = 1$. In particular, the local index $i_{k, 1}$ of $A$ is the integer which is 
uniquely determined by 
\[
i_{k, 1} \in [0, r_{k}-1]\,~\mathit{and}~\, -i_{k, 1}f \equiv 1\,~\mathit{mod}~\, r_{k}\, .
\]
\end{enumerate}
\end{lemma}

\begin{proof}[Proof of $(1)$] By definition and by the fact that $K_{X}$ is a generator of the local class group at each singular point of $X$ (cited above), 
the Gorenstein index $r$ is the least common multiple of the 
orders of the local class groups. Now the assertion (1) follows from the local version cited above. 
\end{proof}

\begin{proof}[Proof of $(2)$] Let $Nr$ be any large multiple of $r$. Since $X$ has only isolated singularities and $NrK_x$ is very ample, by using Bertini's theorem we can find 
a smooth element $S$ in $\vert NrA\vert$. Since $A$ is $\Q$-Cartier and integral, $A \vert S$ is an well-defined integral Cartier divisor on a smooth $S$ 
and 
\[
NrA^{3} = (A.A.S) = (A\vert S)_{S}^{2} \in \Z \, .
\]
Replacing $N$ by $N+1$, we also have that $(N+1)rA^3 \in \Z$. 
Hence $rA^3 \in \Z$. 
The positivity of $n$ follows from the ampleness of $A$. 
\end{proof}

\begin{proof}[Proof of $(3)$] Let $P$ be a singular point of $X$ with local Gorenstein index $r_{P}$. 
For the same reason as before, we can write 
$A = i(P, A)K_{X}$ in the local class group $G_{P}$ at $P$. Since $-K_{X} = fA$, 
we then have $K_{X} = -i(P, A)fK_{X}$ in $G_{P}$. Since $G_{P} = \langle K_{X} \rangle \simeq \Z /r_{P}$, we have $i(P, A)f \equiv 1\, \mathrm{mod}\, r_{P}$. 
Thus, $(f, r_{P}) = 1$ at each $P$. 
Therefore 
$(f, r) = 1$ and $(f, r_{k}) = 1$ by $r = \lcm\, (r_{P})_{P \in 
\mathrm{Sing}\, X} = \lcm\, (r_{k})_{k=1}^{m}$. Since $-K_{X} = fA$ 
at $P_{k}$, we then have that 
\[
-K_{X} = fi_{k, 1}K_{X}\, ,\,\,~\mathrm{i.e.}~\,\,  fi_{k, 1} \equiv 
-1~\mathrm{in}~\, \Z /r_{k}\, .
\] 
This implies the last assertion. 
\end{proof}

\begin{rem}
For each $n$, the local indices $i_{k, n}$ of $nA$ satisfy two obvious relations 
\[
i_{k, n} \equiv ni_{k, 1}\, 
~\mathrm{mod}~\, r_{k}\,~\mathrm{and}~\, 0 \leq i_{k, n} \leq r_{k} - 1\, .
\]
These two relations (together with necessary division and subtraction) 
determine the value $i_{k, n}$ from $i_{k, 1}$. 
\end{rem} 
 
Now we can write down the singular Riemann-Roch formula for 
$\Q$-Fano 3-folds: 

\begin{thm}[\cite{Ka1},~\cite{KMM},~\cite{young}]
Let $X$ be a $\Q$-Fano $3$-fold of Fano index $f = f(X)$ and $A$ be 
a primitive Weil divisor $($so that $-K_{X} = fA$$)$. Let $\mathcal{B}(X) = 
\{P_{k} = [r_{k}, a_{k}]\}_{k=1}^{m}$ be the basket of singularities of $X$. 
Put $P_{n}(X) := \chi(nA)$. We define the Hilbert series of $X$ as the formal 
power series $P(X, t) = \sum_{n = 0}^{\infty} P_{n}(X)t^{n}$.
Then 
\begin{enumerate}
\item[(1)]
\[\chi(\mathcal{O}_{X}) = 1 = \frac{-K_{X}.c_{2}(X)}{24} + \sum_{k=1}^{m}\frac{r_{k}^{2} -1}{24r_{k}}\, .\] 
\item[(2)]  
$P_{n}(X)  = h^{0}(nA)$ for all $n > -f$ and $P_{n}(X) = 0$ for all 
$-f < n < 0$. 
\item[(3)] 

\begin{eqnarray*}\hspace*{1cm}
\chi(\mathcal{O}_{X}(nA)) =  \chi(\mathcal{O}_{X}) + \frac{n(n+f)(2n+f)}{12}A^{3}\hspace*{4.6cm}\\
 + \frac{nA.c_{2}(X)}{12} 
+  \sum_{k=1}^{m} \left( -i_{k, n}\frac{r_{k}^{2} -1}{12r_{k}} 
+ \sum_{j=1}^{i_{k, n}-1}\frac{\overline{b_{k}j}(r_{k}-\overline{b_{k}j}) }{2r_{k}}\right), 
\end{eqnarray*}
where $i_{k, n} \in [0, r_{k}-1]$ is the local index of $nA$ at $P_{k}$, $b_{k} \in [0, r_{k}-1]$ is the integer such that 
$a_{k}b_{k} \equiv 1~\mathit{mod}~r_{k}$ and 
$\overline{b_{k}j} \in [0, r_{k}-1]$ is the integer such that 
$\overline{b_{k}j} \equiv b_{k}j~\mathit{mod}~r_{k}$. In particular,

\begin{eqnarray*}
P(X, t)
& {}=&\frac{1}{1-t}+\frac{(f^{2}+3f+2)t+(-2f^{2}+ 8)t^{2}+(f^{2}-3f + 2)t^{3}}{12(1-t)^{4}}A^{3}\\
& & + \frac{t}{(1-t)^{2}}\frac{A.c_{2}(X)}{12}\\
& & + \sum_{k=1}^{m}
\frac{1}{1-t^{r_{k}}}\left( \sum_{l=1}^{r_{k}-1}\left(-i_{k, l}\frac{r_{k}^{2}-1}{12r_{k}}+
\sum_{j=1}^{i_{k,l}-1}\frac{\overline{b_{k}j}(r_{k}-\overline{b_{k}j})}{2r_{k}}\right)t^{l}\right).
\end{eqnarray*}

\end{enumerate}
\end{thm}
  
\begin{proof} The statement (1) is proved in [Ka1, Section~2] 
(See also [Re, Corollary~10.3]). Since $A$ is ample and 
$nA = K_{X} + (n+f)A$, 
by using 
Kawamata-Viehweg vanishing theorem (see e.g. [KMM, Theorem~1.2.5]), we have \\
$ h^{i} ( \mathcal{O}_{X} ( nA ) ) = 0$ for $ n > -f $ and for $i > 0$. This implies the 
first equality of (2). For the second equality, we may now note that 
$h^{0}(nA) = 0$ for $n < 0$. 
The first equality of (3) is the so-called singular Riemann-Roch formula. 
This is shown by [Re, Theorem~10.2] for an arbitrary projective terminal 
3-fold. Observe that 
$$\frac{t^{N}}{(1-t)^{N+1}} = (-1)^{N}\sum_{L \geq N}\frac{L!}{N!(L-N)!}t^{L}$$ and that $i_{k, n} = i_{k, n + r_{k}}$. Here the last equality 
is a direct consequence of the definition of the index. Now 
the second equality in (3) follows from the first equalities.    
\end{proof}

\begin{cor}
Under the same notation as in Theorem~$1.4$, if $f = f(X) \geq 3$, then 
\begin{eqnarray*}
\lefteqn{A^{3}= \frac{12}{(f-1)(f-2)}}\hspace*{12.5cm}\\
{}\times \left(1 - \frac{A.c_{2}(X)}{12} + \sum_{k=1}^{m} \left(-i_{k, -1}\frac{r_{k}^{2} -1}{12r_{k}} +  
\sum_{j=1}^{i_{k, -1}-1}\frac{\overline{b_{k}j}(r_{k}-\overline{b_{k}j})}{2r_{k}}\right)\right)\, .
\end{eqnarray*}
\end{cor}

\begin{proof} Since $f \geq 3$, we have $\chi(\Oh_{X}(-A)) = 0$ by the theorem~1.4~(2). Now, 
substituting $n = -1$ into the first equality of the theorem~1.4~(3), we get the result. 
\end{proof}

Next, we recall the boundedness theorem due to Kawamata. In his paper [Ka3], he shows the following: 

\begin{thm} [\textrm{[Ka3, ~Proposition~1 (see also Theorem~2)]}]
There is a universal constant $b > 0$ 
such that 
$$(-K_{X})^{3} \leq b(-K_{X}.c_{2}(X))$$ 
holds for all $\Q$-Fano $3$-folds $X$. In particular, 
$0 < \left(-K_{X}.c_{2}(X)\right)$. 
\end{thm}

However, he shows more in the course of proof, as we shall now explain. 
Let $X$ be a $\Q$-Fano 3-fold and  
$\E := (\Omega_{X}^{1})^{**}$ be 
the double dual of the sheaf of K\"ahler differentials of $X$. 
If $\E$ is not $\mu$-semistable (with respect to $-K_{X}$), 
then one can take the so-called maximal destabilizing sheaf  $\F$ of $\E$, i.e. a (unique) $\mu$-semistable subsheaf $\F \subset \E$, 
which is necessarily reflexive and of rank $s = 1$ or $2$, such that 
\[
\frac{(c_{1}(\F).(-K_{X})^{2})}{s} > \frac{(K_{X}.(-K_{X})^{2})}{3}\, .
\]

Set $c_{1}(\F) = tK_{X}$. 
For this expression, we used the fact that 
$\rho(X) = 1$ and $X$ is $\Q$-factorial. It is shown that 
$0 < t < s/3$ [Ka3, Pages 442-443]. 
Under these notations, one can say what he showed 
as in the following more effective form: 

\begin{thm}[\textrm{[Ka3,~Proposition $1$ (see also the proof there)]}]
Under the 
above setting, one has:

\begin{enumerate}  
\item[(1)] If $\E$ is $\mu$-semistable, then 
\[(-K_{X})^{3} \leq 3(-K_{X}.c_{2}(X))\, .\] 

\item[(2)]  If $\E$ is not $\mu$-semistable and $s = 1$, then one of the following holds:
\[(1-t)(1+3t)(-K_{X})^{3} \leq 4(-K_{X}.c_{2}(X))\, ,\] 
or 
\[(tu+(t+u)(1-t-u))(-K_{X})^{3} \leq (-K_{X}.c_{2}(X))\] 
for some rational number $u$ such that $t < u < 1-t-u$.

\item[(3)] If $\E$ is not $\mu$-semistable and $s = 2$, 
then  
\[t(4-3t)(-K_{X})^{3} \leq 4(-K_{X}.c_{2}(X))\, .\] 
\end{enumerate} 
\end{thm} 

\section{ Fano indices of $\Q$-Fano 3-folds}

In this section, we shall show Theorem 0.3. Throughout this section, we assume 
that $X$ is a $\Q$-Fano 3-fold of Fano index $f = f(X) \geq 3$, 
$\mathcal{B}(X) = \{[r_{k}, a_{k}]\}_{k=1}^{m}$ is 
the basket of singularities of $X$, 
$r := r(X) (= \lcm(r_{k})_{k=1}^{m})$ 
is the Gorenstein index and $A$ be a primitive Weil divisor. 
\par 
\vskip 4pt 

The following quantity is important in the sequel:

\begin{definition} 
\[ 
B(X) := 
r\left(24 - \sum_{k=1}^{m} \left(r_{k} - \frac{1}{r_{k}}\right)\right)\,. 
\] 
\end{definition}

Using Theorem 1.7, we shall first deduce the following inequality:

\begin{prop} 
\[ 
(4f^{2} - 3f)A^{3} \leq 4(-K_{X}.c_{2}(X))\, .
\]
\end{prop} 

\begin{proof} We show the inequality by dividing into the four cases 
in Theorem 1.7. 
\par 
\textit{ $1$. the case where} $(-K_{X})^{3} \leq 3\left(-K_{X}.c_{2}(X)\right)$. 
\par
In this case we have
\[\frac{4f^{3}}{3}A^{3} \leq 4(-K_{X}.c_{2}(X))\]
by $-K_{X} = fA$. Since $f \geq 3$, one has also
$$4f^{2} - 3f \leq 4f^{2} \leq \frac{4f^{3}}{3}\, .$$
Combining these two inequalities, we get the desired inequality. 
\par 
\textit{ $2$. the case where} $(1-t)(1+3t)(-K_{X})^{3} \leq 4(-K_{X}.c_{2}(X))$. 
\par 
In this case we have
\[(1-t)(1+3t)f^{3}A^{3} \leq 4(-K_{X}.c_{2}(X))\, .\]
by $-K_{X} = fA$. 
Since the function $(1-x)(1+3x)$ is increasing in the range 
$0 \leq x \leq 1/3$ and since $0 < t <1/3$, we have 
\[f^{3}A^{3} \leq 4(-K_{X}.c_{2}(X))\, .\]
One has also
$$4f^{2} - 3f \leq 4f^{2} \leq f^{3}~(f \geq 4)\,\, \text{and}\,\,  
4f^{2} - 3f = f^{3}~(f = 3)\, .$$
Thus, we get the desired inequality in this case, too. 
\par 
\textit{$3$. the case where} $(tu+(t+u)(1-t-u))(-K_{X})^{3} 
\leq (-K_{X}.c_{2}(X))$ 
\textit{ for some rational number} $u$ {\it such that} $t < u < 1-t-u$. 
\par
By $0 < t < u < 1-t-u$, we have also $0 < t < 1-2t$, i.e. $0 < t < 1/3$ 
and $t < u < (1-t)/2$. 
Since 
\[
tu+(t+u)(1-t-u) = -\left(u - \frac{1-t}{2}\right)^{2} + \frac{-3t^{2} + 2t + 1}{4}\, ,
\]
$tu +(t+u)(1-t-u)$ is increasing with respect to $u$ in the range 
$(t, (1-t)/2)$. 
Thus 
$$-3t^{2} + 2t = t^{2} + 2t(1-2t) \leq tu+(t+u)(1-t-u)\, ,$$ 
and therefore  
\[
(-12t^{2} + 8t)f^{3}A^{3} \leq 4(-K_{X}.c_{2}(X))\, .\]
Since 
\[
t \in \left\{ \frac{1}f, \frac{2}f, \dots , \frac{n}f, \dots \right\} \cap 
\left( 0, \frac13\right) \, ,
\] 
we have 
\[
-12t^{2} + 8t\, = -12\left(t-\frac{1}{3}\right)^{2} + \frac{4}{3}\, \geq -\frac{12}{f^{2}} + \frac{8}{f}\, .
\]
Therefore 
\[
(-12t^{2} + 8t)f^{3} \geq 8f^{2} -12f \geq 4f^{2} - 3f\]

for $f \geq 3$. From this inequality, we obtain 
\[
(4f^{2} - 3f)A^{3} \leq 4(-K_{X}.c_{2}(X)\, ).\]

\par

\textit{ $4$. the case where} $t(4-3t)(-K_{X})^{3} 
\leq 4(-K_{X}.c_{2}(X))$.
\par
By $-K_{X} = fA$, we have 
\[t(4-3t)f^{3}A^{3} \leq 4(-K_{X}.c_{2}(X))\, .\] 
Again by $-K_{X} = fA$ and by $c_{1}(\F) = tK_{X}$ (by the definition of $t$), 
one has $c_{1}(\F) = -ftA$ in the Weil divisor class group (in the numerical sense). Since $A$ is 
a generator of this group (by the $\Q$-factoriality of $X$ and $\rho(X) = 1$), we have $ft \in \Z$. Therefore 
\[
t \in \left\{\frac{1}{f}, \frac{2}{f}, \cdots , \frac{n}{f}, \cdots \right\} 
\cap \left(0, \frac{2}{3}\right)\, .
\]  
Using this, we obtain 
\[
t(4-3t)\, = -3\left(t - \frac{1}{3}\right)^{2} + \frac{1}{3}\, \geq \frac{1}{f}
\left(4 - \frac{3}{f}\right)\, .
\] 
Therefore 
\[
(4f^{2} -3f)A^{3} = \frac{1}{f}\left(4 - \frac{3}{f}\right)f^{3}A^{3} 
\leq t(4-3t)f^{3}A^{3} \leq 4(-K_{X}.c_{2}(X))\, .
\]
 
Now we are done. 
\end{proof}

The next inequality is crucial for us. 

\begin{cor}
\[ 
4f^{2} - 3f \leq 4B(X)\, .
\]
\end{cor}

\begin{proof} We have $A^{3} \geq 1/r$ by Lemma 1.2~(2). 
Substituting this inequality and the equality in Theorem 1.4~(1) 
into the inequality of Proposition~2.2, 
we obtain the result.  
\end{proof} 

First we shall bound $B(X)$ from the above, then one can also 
estimate $f$. 

\begin{prop}

\begin{enumerate}
\item[(1)]$($\textup{\cite[Proof~of~Theorem~2]{kawamata}, \cite[Proof~of~1.2 (1)]{KMMT}}$)$
\[ 
\sum_{k=1}^{m}\left(r_{k} - \frac{1}{r_{k}}\right) < 24\, .
\]
\item[(2)] 
\[
 0 < B(X) \leq 2489\, .
 \]
Moreover, the right equality in $(2)$ holds if and only if $($$m = 4$ and~$)$ 
\[ 
\{r_{1}, r_{2}, r_{3}, r_{4}\} = \{3, 4, 5, 7\}\, 
\]
\end{enumerate} 

\end{prop} 

\begin{rem}
For a $\Q$-Fano 3-fold $\PP (3,4,5,7)$, we have 
\[
\{r_{k}\}_{k=1}^{m} = \{3, 4, 5, 7\}\,\, \text{and}\,\, B(\PP (3,4,5,7)) 
= 2489\, .
\]
Note also that $f(\PP (3,4,5,7)) = 19$. ($\mathrm{c.f.}$~Proposition~2.13.)
This already indicates that the value $f = 19$ is something special.  
\end{rem} 

\begin{proof} Our argument here is suggested by T. Katsura.
\par  
Since $(-K_{X}.c_{2}(X)) > 0$ by Theorem 1.6 and since $\chi(\Oh_{X}) = 1$, 
we have the first inequality 
\[ 
\sum_{k=1}^{m}\left(r_{k} - \frac{1}{r_{k}}\right) < 24\,
\] 
by Theorem~1.4(1). This is also equivalent to $0 < B(X)$. In what follows, 
we seek the maximum value of $B(X)$ together 
with the basket which attains the maximum. 
For this purpose, 
it is more convenient to observe first the following purely arithmetical claim (apart from $\Q$-Fano 3-folds for a moment):

\begin{lemma} Let $\{r_{k}\}_{k=1}^{m}$ be a finite sequence 
of integers such that $r_{k} \geq 2$ for all $k$ and such that 
$$\sum_{k=1}^{m}\left(r_{k} - \frac{1}{r_{k}}\right) < 24\, .$$ 
Set $r = \lcm\, (r_{k})_{k=1}^{m}$. Then 
$$B(\{r_{k}\}_{k=1}^{m}) := r\left(24 - \sum_{k=1}^{m}\left(r_{k} - \frac{1}{r_{k}}\right)\right) \leq 2489$$ 
and the equality holds if and only if $(m = 4$ and~$)$ 
$\{r_{1}, r_{2}, r_{3}, r_{4}\} = \{3, 4, 5, 7\}$.
\end{lemma}

\begin{proof}[Proof of Lemma]
By the second condition of $\{f_k\}_{r=1}^m$, 
we have $r_{k} \leq 24$ for all $k$ and 
$$m\cdot\frac{3}{2} \leq \sum_{k=1}^{m}\left(r_{k} - \frac{1}{r_{k}}\right) < 24\, .$$ 
Thus $m \leq 15$. In particular, there are only finitely many sequences 
$\{r_{k}\}_{k=1}^{m}$ which satisfy the initial two conditions. 
So, there is certainly the maximum of $B(\{r_{k}\}_{k=1}^{m})$, say 
$M$, when $\{r_{k}\}_{k=1}^{m}$ varies. In what follows, we seek the value $M$ as well as the sequences which attain the maximum. 

\begin{cl}The sequence $\{3, 4, 5, 7\}$ satisfies the initial conditions and \\
$B( \{3, 4, 5, 7\} ) = 2489$. In particular, 
$2489 \leq M$. 
\end{cl} 

\begin{proof} This follows from a direct calculation. 
\end{proof}

\begin{cl}If $r_{j} = p^{a}q^{b}$ $(p$ and $q$ are different prime numbers, $a \geq 1$ and $b \geq 1)$ for some $j$, then 
$B(\{r_{k}\}_{k=1}^{m}) < M$. 
\end{cl}
 
\begin{proof} We may assume that $j = m$, i.e. $r_{m} = p^{a}q^{b}$. Consider a new sequence $\{s_{k}\}_{k=1}^{m+1}$ defined by 
\[
s_{k} = r_{k}\,~\textrm{for}\,~k \leq m-1\,~\textrm{and}~s_{m} = p^{a}, s_{m+1} = q^{b}\, .
\] 
Then $\lcm\,(s_{k})_{k=1}^{m+1} = \lcm\,(r_{k})_{k=1}^{m}$ and 
\begin{eqnarray*}
\lefteqn{\left(p^{a}q^{b} - \frac{1}{p^{a}q^{b}}\right) - 
\left(\left(p^{a} - \frac{1}{p^{a}}\right) + \left(q^{b} - \frac{1}{q^{b}}\right)\right)}\hspace*{5.7cm}\\
= (p^{a} -1)(q^{b} - 1)\left(1 -\frac{1}{p^{a}q^{b}}\right) > 0\, .
\end{eqnarray*}
Thus, the sequence $\{s_{k}\}_{k=1}^{m+1}$ satisfies the initial conditions and 
\[
B(\{r_{k}\}_{k=1}^{m}) < B(\{s_{k}\}_{k=1}^{m+1})\, .\,\, 
\]
\end{proof}

\begin{cl} If there are two numbers $i \not= j$ such that 
$r_{i} = p^{a}$ and $r_{j} = p^{b}$ $(p$ is a prime number and $a \geq b$ 
are positive integers $)$, then 
$B(\{r_{k}\}_{k=1}^{m}) < M$. 
\end{cl}
 
\begin{proof} As before, we may assume that $r_{m-1} = p^{a}$ and $r_{m} = p^{b}$. 
Consider a new sequence $\{s_{k}\}_{k=1}^{m-1}$ defined by 
$$s_{k} = r_{k}\, .$$ 
Then $\lcm\,(s_{k})_{k=1}^{m-1} = \lcm\,(r_{k})_{k=1}^{m}$ and 
\[
\sum_{k=1}^{m-1} \left(s_{k} - \frac{1}{s_{k}}\right) < \sum_{k=1}^{m} \left(r_{k} - 
\frac{1}{r_{k}}\right).\]
 
Thus, the sequence $\{s_{k}\}_{k=1}^{m-1}$ satisfies the condition and 
\[
B(\{r_{k}\}_{k=1}^{m}) < B(\{s_{k}\}_{k=1}^{m-1})\, .\,\, 
\] 
\end{proof}
 
By Claims~2.8 and~2.9, we may now assume that all $r_{i}$ are primary, i.e. 
$r_{i} = p_{i}^{a_{i}}$ where $p_{i}$ is a prime number, and $r_{i}$ are coprime to one another, i.e. $p_{i} \not= p_{j}$ if $i \not= j$. In particular, 
$\lcm\, (r_{k})_{k=1}^{m} = \prod_{k=1}^{m}r_{k}$. 

\begin{cl} If $m \not= 4$, then $B\left( \{r_{k}\}_{k=1}^{m}\right) < M$. 
\end{cl} 

\begin{proof} If $m \geq 5$, then by $r_{k} \leq 24$ and by 
the coprime conditions above (now assumed), we have 
\begin{eqnarray*}
& &\sum_{k=1}^{m} \left(r_{k} - \frac{1}{r_{k}}\right) 
\\
& &\geq 
\left(2 - \frac{1}{2}\right) + \left(3 - \frac{1}{3}\right) + 
\left(5 - \frac{1}{5}\right) + \left(7 - \frac{1}{7}\right)
+ \left(11 - \frac{1}{11}\right)  \hspace*{-1cm}\\
& &\geq 2 + 3 + 5 + 7 + 11 - \frac{1}{2}\cdot 5 > 24\, ,  
\end{eqnarray*}
a contradiction to the initial conditions. Therefore $m \leq 4$. 
\par \vskip 4pt
Next assume that $m \leq 3$. Since the value $B(\{r_{k}\}_{k=1}^{m})$ is invariant even if we formally 
add terms $1$ into the sequence, we may assume that the sequence 
is of the form 
\[ 
\{r_{k}\}_{k=1}^{m} = \{a, b, c\} 
\]
in which some of $a$, $b$, $c$ are allowed to be $1$. 
By the coprime condition, we have $\lcm(r_{k})_{k=1}^{m} = abc$ and 
\begin{eqnarray*}
B(\{r_{k}\}_{k=1}^{m})& =& abc\left(24 - a - b - c + \frac{1}{a} + 
\frac{1}{b} + \frac{1}{c}\right) \\
&\leq& abc(27 - a - b - c) \leq abc(28 - a - b- c) \\
&\leq &\left(\frac{28}{4}\right)^{4} = 2401 < 2489 \leq M\, .
\end{eqnarray*}
Hence the claim follows. 
\end{proof}

Now we may furthermore assume that $m = 4$ and (without loss of generality) 
that 
$$r_{1} > r_{2} > r_{3} > r_{4}\, .$$  

\begin{cl} 
If $r_{1} \geq 11$, then $B(\{r_{k}\}_{k=1}^{4}) < M$. 
\end{cl}

\begin{proof} By coprime condition and $r_{1} \geq11$, we have 
\begin{eqnarray*}
B(\{r_{k}\}_{k=1}^{4})\hspace*{-5pt}&=&r_{1}r_{2}r_{3}r_{4}\left(24 - r_{1} - r_{2} - r_{3} - r_{4} + \frac{1}{r_{1}} + 
\frac{1}{r_{2}} + \frac{1}{r_{3}} + \frac{1}{r_{4}}\right)\\
&\leq& r_{1}r_{2}r_{3}r_{4}((28 - r_{1}) - r_{2} - r_{3} - r_{4})\\ 
&\leq& r_{1}\left(\frac{28 - r_{1}}{3}\right)^{3} \leq 11 \cdot 6^{3} = 2376 < 
2489 \leq M\, .
\end{eqnarray*}
This implies the claim. 
\end{proof}

Now we can complete the proof of Lemma~2.6. By Claim~2.11 and coprime conditions, we have $r_{1} \leq 9$ for the maximum $B$. 
Now there are exactly $6$ sequences $r_{1} > r_{2} > r_{3} > r_{4}$ 
which satisfy the coprime conditions and $r_{1} \leq 9$. They are: 
$$(9, 8, 7, 5)\, ,\, (9, 7, 5, 4)\, ,\, (9, 7, 5, 2)\,\,$$ 
$$(8, 7, 5, 3)\, ,\, (7, 5, 4, 3)\, ,\, (7, 5, 3, 2).$$ 
Among these six candidates, the first two sequences do not satisfy the initial condition 
$\sum_{k=1}^{m}(r_{k} - \frac{1}{r_{k}}) < 24$. Now, by calculating $B(\{r_{k}\}_{k=1}^{4})$ for the other four sequences, 
we obtain the desired result as in Lemma~2.6. \end{proof}
Now Proposition 2.4 follows from Lemma~2.6. \end{proof}
\par 
\vskip 4pt
By combining Corollary~2.3 and Proposition~2.4 (2), we obtain a rough estimate of $f$: 

\begin{cor} $f \leq 50$.  
\end{cor}
 
\begin{proof} By Corollary~2.3 and Proposition 2.4 (2), we have 
$$4f^{2} - 3f \leq 4 \cdot 2489 = 9956\, .$$ 
This implies $f \leq 50$. 
\end{proof}
In order to obtain an optimal estimate $f \leq 19$, we need one more work. What we will do from now is to seek integral solutions of a system of equalities and inequalities which the baskets $B(X)=\{[r_k, a_k]\}^m_{k=1}$ of $\Q$-Fano 
$3$-folds $X$ must satisfy. 
\par 
\vskip 4pt 
If there is a $\Q$-Fano 3-fold $X$ of Fano index $f~(\geq 3)$, 
for which we now know $f \le 50$, 
there must be integer solutions $m$, $r_{k}$, $a_{k}$ (or equivalently $b_{k}$), $i_{k, n}$ 
($1 \leq k \leq m$)  of the following equations and inequalities: 
\begin{enumerate} 
\item[(1)] By Proposition~2.4~(2) 
\[ 
\sum_{k=1}^{m}\left(r_{k} - \frac{1}{r_{k}}\right) < 24\,~\textrm{and}~(r_k, a_k)=1.
\] 
\item[(2)]  By Lemma~1.2~(3), for all $k$, 
\[ (f, r_{k}) = 1\, .\]
\item[(3)] By Corollary~1.5,
\begin{eqnarray*} 
& &A^{3}
=  \frac{12}{(f-1)(f-2)}\\
& &\phantom{\sum}\ \times\left(1-\frac{A.c_{2}(X)}{12} + \sum_{k=1}^{m} \left(-i_{k,-1}\frac{r_{k}^{2} -1}{12r_{k}} + \sum_{j=1}^{i_{k,-1}-1}\frac{\overline{b_{k}j}(r_{k}-\overline{b_{k}j})}{2r_{k}}\right)\right)\\
& &\phantom{A^3}>0\,.
\end{eqnarray*}

\item[(4)] By Theorem~1.4~(2), for all $n = -1, -2, \cdots, -(f-1)$
\begin{eqnarray*}
\chi(\Oh_{X}(nA)) \hspace*{8cm}\\= 1 + \frac{n(n+f)(2n+f)}{12}A^{3}
+ \frac{nA.c_{2}(X)}{12} \hspace*{3.8cm}\\
{} \,\,\,\,\,\,
+ \sum_{k=1}^{m}\left(-i_{k, n}\frac{r_{k}^{2} -1}{12r_{k}} + \sum_{j=1}^{i_{k, n}-1}\frac{\overline{b_{k}j}(r_{k}-\overline{b_{k}j})}{2r_{k}}\right) \\= 0\, .
\end{eqnarray*} 

\item[(5)] By Proposition~2.2
\[ (4f^{2} - 3f)A^3 \leq 4(-K_X.c_{2}(X)
)\, .\]
\end{enumerate}
 
Here by 1.4(1), we have 
$$-K_X.c_{2}(X) =24-\sum_{k=1}^{m}\left(r_{k} - \frac{1}{r_{k}}\right)\, .$$ 
As we remarked before, there are 
only finitely many integers  $\{r_{k}\}_{k=1}^{m}$ which satisfy the inequality (1). 
For each $f$ and $\{r_{k}\}_{k=1}^{m}$, 
there is a unique integer of $i_{k, 1}$ (whence 
$i_{k, n}$) by Lemma~1.2~(3) and Theorem~1.4, and finitely many integers $b_{k}$ (or equivalently $a_{k}$) by $b_{k} \in [0, r_{k}-1]$. 
For each such possibility in the range 
$9 \leq f \leq 50$, we check if it satisfies (2)--(5), by additions and multiplications. 
In principle, we can do this by hand. However, it is a little messy to do so and we use a computer program Magma. 
Among 5 conditions, the condition (4) seems fairly strong. 
As a result, we actually find that there are no integer solutions when $f \geq 20$. 
Thus we have $f \leq 19$. 

For instance, our programs give the following list(c.f. Tabel~1) : 
\vskip 4pt

\begin{table}[h]
 \[
 \renewcommand{\arraystretch}{1.5}
 \begin{array}{|c|c|c|c|c|c|}
\hline
  f & \sharp~(1) & \sharp~(2) & \sharp~(3) & \sharp~(4) & \sharp~(5)  \\
\hline\hline
13 & 25161 & 23187 & 6622 & 6 & 2 \\
\hline
19 & 25161  & 24972 & 7173 & 1  & 1 \\
\hline
20 & 25161 & 714 & 417 & 0 & 0 \\
\hline
23 & 25161 & 25139 & 9261 &  0 & 0 \\
\hline
24 & 25161 & 478 & 329 & 0 & 0 \\
\hline
50 & 25161 & 714 & 167 & 0 & 0\\
\hline
\end{array}
\]
\caption{The number of decreasing}
\end{table}
\vskip 4pt

In order to make this process clear, 
we gave programs we used in the appendix. 

The next proposition shows the optimality of the estimate $f \leq19$: 

\begin{prop} 
The weighted projective space 
$\PP (3,4,5,7)$ is a $\Q$-Fano $3$-fold of Fano index $19$.  
\end{prop} 

\begin{proof} Recall that $X := \PP (3, 4, 5, 7)$ 
is an abelian quotient of $\PP^{4}$ by an obvious action by the abelian group 
$C_{3} \times C_{4} \times C_{5} \times C_{7}$. 
Therefore, $X$ is $\Q$-factorial, 
the Weil divisor class group is generated by the Serre's twisting sheaf 
$\Oh_{X}(1)$, and 
that $\text{Sing}(X) = \{[3,1], [4,1], [5, 2], [7,2]\}$, which are terminal.  
In addition, by the canonical bundle formula (See [Do]), we have 
\[ 
\Oh_{X}(K_{X}) = \Oh_{X}(-19)\, .
\]
Thus the Fano index of $X$ is $19$. 
\end{proof}

However, since we obtain all the solutions of (1)--(5) for 
$3 \le f \le 50$, we can say more about $X$ for each possible $3 \leq f \leq 19$. For instance,  
We also find that if $f = 19$, then $X$ necessarily satisfies
\[
A^{3} = \frac{1}{3\cdot 4 \cdot 5 \cdot 7}\, ,
\]
\[
\mathcal{B}(X) = \{[3,1], [4,1], [5, 2], [7,2]\}\, ,
\] 
and 
\[ 
P(X, t) = \frac{1}{(1-t^{3})(1-t^{4})(1-t^{5})(1-t^{7})}\, .
\] 

In this way, we can obtain the assertion (1) and (2) in the Theorem 0.3. 
The assertion (3) is now easily proved. Let $f$ be an integer in the assertion (3). 
Then one can actually construct explicit examples of $\Q$-Fano 3-folds of index $f$ 
as general hypersurfaces in weighted projective spaces. For instance, we have the following simple examples (in which the equations are chosen to be general) with indicated singular points:

\begin{table}[h]
 \[
 \renewcommand{\arraystretch}{1.5}
 \begin{array}{|c|c|c|}
\hline
  f & X & \mathrm{Sing}(X)=\mathcal{B}(X)  \\
\hline\hline
17 &  \PP (2,3,5,7) & \{[2,1], [3,1], [5,1], [7,3]\} \\ 
\hline
13 & \PP(1,3,4,5) & \{[3,1], [4,1], [5,2]\} \\  
\hline
11 & \PP(1, 2,3,5) & \{[2,1], [3,1], [5,2]\} \\
\hline
9  & (6) \subset \PP (1,2,3,4,5) & \{[2,1], [4,1], [5,2]\}\\ 
\hline
8  &  (6) \subset \PP (1,2,3,3,5) & \{[3,1], [3,1], [5,1]\} \\
\hline 
7  & \PP (1,1,2,3) &  \{[2,1], [3,1]\} \\  
\hline
 6 & (6) \subset \PP (1, 1, 2,3,5) & \{[5,2]\} \\ 
\hline
5  & \PP (1,1,1,2) & \{[2,1]\} \\
\hline
4  & \PP^{3} &  \emptyset \\ 
\hline
3 & (2) \subset \PP^{4} & \emptyset \\
\hline
2 & (3) \subset \PP^{4} & \emptyset \\
\hline 
1 & (4) \subset \PP^{5} & \emptyset \\ 
\hline
\end{array}
\]
\caption{Examples of $\Q$-Fano~3-folds}
\end{table}
\vskip 4pt
Now we are done. Q.E.D. for the Theorem~0.3.  \qed

\begin{rem} Similarly, using Magma program, 
(but use $\chi(A) \geq 0$ instead of $A^{3} > 0$ when $f = 1,2$),   we obtain
$$-K_{X}^{3} \leq \frac{2 \cdot 5^{3}}{3}$$ 
for a $\Q$-Fano 3-folds $X$. Unfortunately, we do not know whether this estimate is optimal or not.
\end{rem}

Finally, we pose three interesting unsettled problems which are 
closely related to our theorem: 
\par \vskip 4pt

{\em Question~$1$. }~Is there a more intrinsic reason why 
$f(X) \leq 19$ should hold?  
\par \vskip 4pt

{\em Question~$2$. }~Classify all $X$ with $f(X) = 19$ 
up to isomorphism. \\
\hspace*{66pt}
$X \simeq \PP (3,4,5,7)$ if $f(X) = 19$? 
\par \vskip 4pt

{\em Question~$3$. }~Is there a $\Q$-Fano 3-fold of Fano index $10$? 

\section{Appendix: Magma Program} 
This is  a program which we used at the final step of the proof 
for the calculation of 
$1/12Ac_2(X)$, $A^3$ and $\mathit{P_n(X)}$ and the Hilbert series $P(X,t)$ 
from Fano index, local indices, and the baskets of singularities, 
i.e. from the values $f$, $i_{n,k}$, $[r_k, a_k]$. 

Here, $\mathrm{BB}$ is a list of 
Baskets generated automatically by computer under the condition of 
$\sum^m_{k=1} \left(r_k-1/r_k \right)<24$.

\vspace{1cm}
\begin{verbatim}

//////////////////////////////////////////////////////
// Build the Hilbert series
//////////////////////////////////////////////////////

forward Ac2over12_is, contribution;
intrinsic FanoHilbertSeries(f::RngIntElt,B::SeqEnum) 
-> RngElt
{The Hilbert series of a Fano 3-fold of Fano index f
 and basket B}

    K := RationalFunctionField(Rationals());
    t := K.1;
    I := 1/(1-t);
    II := 1/12*A3_is(f,B)*
        ((f^2+3*f+2)*t+(-2*f^2+8)*t^2+(f^2-3*f+2)*t^3) 
        /(1-t)^4;
    III := Ac2over12_is(f,B)*t/(1-t)^2;
    IV := &+[ Parent(t) |
        &+[ Parent(t) | contribution(f,r,a,n)*t^n 
         	    : n in [1..r-1] ] / (1-t^r)
                where r is p[1]
                where a is p[2] : p in B ];
    return I + II + III + IV;
end intrinsic;

//////////////////////////////////////////////////////
// Auxiliary functions
//////////////////////////////////////////////////////

function i_is(f,r,n)
    h,u,v := XGCD(f,r);
    return (-n*u) mod r;
end function;

bar := func< m,r | m mod r >;


inv := func< a,r | i_is(a,r,1) >;

function contribution(f,r,a,n)
    i := i_is(f,r,n);
        b := inv(a,r);
    first := -i*(r^2-1)/(12*r);
    if i in {0,1} then
        extra := 0;
    else
        extra := &+[ bar(b*j,r)*(r-bar(b*j,r))/(2*r) :
                         j in [0..i-1] ];
    end if;
    return first + extra;
end function;

function Ac2over12_is(f,B)
    sumpart := &+[ Rationals() | (r^2-1)/(12*r) where r is p[1] 
					 : p in B ];
    return (2-sumpart)/f;
end function;

// require f ge 3: ...
function A3_is(f,B)
       factor := 12/((f-1)*(f-2));
        c2_part := Ac2over12_is(f,B);
        periodic := &+[ Rationals() | contribution(f,r,a,-1)
                                 where a is p[2]
                                 where r is p[1] : p in B ];
        return factor * (1 - c2_part + periodic);
end function;

intrinsic FanoCoefficient(f::RngIntElt,B::SeqEnum,n::RngIntElt) 
-> RngElt
{The n-th coefficient of the Hilbert series of Fano 
                             with Fano index f and basket B}

    V := 1+1/12*A3_is(f,B)*n*(n+f)*(2n+f)+n*Ac2over12_is(f,B)+
         &+[Rationals()| contribution(f,r,a,n) 
         	     	         where a is p[2]
                                 where r is p[1] :  p in B];
    vprintf User1: "\tP_(%o) = %o\n",n,V;
    return V;
end intrinsic;


 BB := Baskets(24);                                             
 B1 := [ B : B in BB | &and[ GCD(p[1],f) eq 1 : p in B ] ];     
 B2 := [ B : B in B1 | A3_is(f,B) gt 0 ];     
 B3 := [ B : B in B2 | A3_is(f,B)*(4*f^2-3) 
 							 le 48*f*Ac2over12_is(f,B)];       
 Bfinal := [ B : B in B3 | &and[ coeff(f,B,n) eq 0 :
                             n in [-(f-1)..-1]]];

\end{verbatim}

\bibliographystyle{amspalpha}

\end{document}